\documentclass[11pt]{article}
\usepackage{amssymb,latexsym,geometry,tikz,setspace,cite}
\newtheorem{theorem}{Theorem}

\newtheorem{lemma}[theorem]{Lemma}

\begin{document}

\onehalfspace

\title{The Cycle Spectrum of Claw-free Hamiltonian Graphs}

\author{Jonas Eckert, Felix Joos, Dieter Rautenbach}

\date{}

\maketitle

\vspace{-1cm}

\begin{center}
Institut f\"{u}r Optimierung und Operations Research, 
Universit\"{a}t Ulm, Ulm, Germany\\
\{\texttt{jonas.eckert, felix.joos, dieter.rautenbach}\}\texttt{@uni-ulm.de}
\end{center}

\begin{abstract}
If $G$ is a claw-free hamiltonian graph of order $n$ and maximum degree $\Delta$ with $\Delta\geq 24$, then
$G$ has cycles of at least $\min\left\{ n,\left\lceil\frac{3}{2}\Delta\right\rceil\right\}-2$ many different lengths.
\end{abstract}

{\small \textbf{Keywords:} Hamiltonian cycle; cycle spectrum; claw-free graph}\\
\indent {\small \textbf{AMS subject classification:} 05C38} % Paths and cycles

\section{Introduction}

We consider finite, simple, and undirected graphs and use standard terminology \cite{ha}.
For a graph $G$, the set of cycle lengths of $G$ is the {\it cycle spectrum} $S(G)$
and $s(G)=|S(G)|$ is the number of different cycle lengths of $G$.
If $G$ has order $n$ and $n\in S(G)$, 
then $G$ is {\it hamiltonian} and a cycle that contains all $n$ vertices of $G$ is a {\it hamiltonian cycle} of $G$.
If $S(G)=\{ 3,\ldots,n\}$, then $G$ is {\it pancyclic},
and,
if $S(G)=\{ 3,\ldots,\max S(G)\}$, then $G$ is {\it subpancyclic}.

Cycles in graphs are among the most established topics in all of graph theory \cite{ha}
and cycles in claw-free graphs have received special attention \cite{faflry}.
Sufficient conditions for hamiltonicity of claw-free graphs \cite{faflry,sh}
or sufficient conditions for hamiltonicity that involve degree conditions 
on induced claws \cite{brrysh,chzhqi} have been proposed. 
The cycle spectrum of sparse claw-free hamiltonian graphs was considered \cite{muraresa}.
Cycle lengths in hamiltonian graphs with few vertices of large degree \cite{koma,fejaha,hasc,mawo,ma}
and degree conditions for (sub)pancyclism of claw-free graphs \cite{gopf,trveve} were studied.

For a hamiltonian graph $G$ of order $n$ and maximum degree $\Delta$,
Marczyk and Wo\'{z}niak \cite{mawo}
showed the best possible inequalities
$$
s(G)\geq 
\left\{
\begin{array}{rl}
\Delta-1, & \mbox{ for }\Delta\leq \frac{n}{2}\mbox{ and}\\
\frac{n}{2}+\frac{\Delta}{2}-\frac{3}{2}, & \mbox{ for }\Delta>\frac{n}{2}.
\end{array}
\right.$$
In the present paper we present the following improvement of these estimates for claw-free graphs.

\begin{theorem}\label{theorem1}
If $G$ is a claw-free hamiltonian graph of order $n$ and maximum degree $\Delta$ with $\Delta\geq 24$, then
$$s(G)\geq \min\left\{ n,\left\lceil\frac{3}{2}\Delta\right\rceil\right\}-2.$$
\end{theorem}
If $n$ and $\Delta$ are integers with $n-2\geq \Delta\geq 2$ 
and the graph $G$ arises 
from the disjoint union of 
a vertex $v_0$, 
a clique $C_1$ of order $\left\lceil\frac{\Delta}{2}\right\rceil$,
a clique $C_2$ of order $\left\lfloor\frac{\Delta}{2}\right\rfloor$, and 
a path $P$ of order $n-\Delta-1$
by adding an edge between $v_0$ and each vertex in $C_1\cup C_2$,
an edge between one endvertex of $P$ and a vertex in $C_1$, and 
an edge between the other endvertex of $P$ and a vertex in $C_2$,
then $G$ is a claw-free hamiltonian graph of order $n$ and maximum degree $\Delta$.
Furthermore,
$C(G)=\left\{ 3,\ldots,\left\lceil\frac{3}{2}\Delta\right\rceil+1\right\}\cup \{ n-\Delta+2,\ldots,n\}$,
that is, Theorem \ref{theorem1} is tight.

The rest of the paper is devoted to the proof of Theorem \ref{theorem1}.

\section{Preliminary Results and Proofs}

Before we proceed to the proof of our result, we collect some auxiliary statements.
We use standard terminology and notation \cite{ha}.

\begin{theorem}[Hakimi and Schmeichel \cite{hasc}]\label{theoremhasc}
Let $G$ be a graph with a hamiltonian cycle $v_0v_1\ldots v_{n-1}v_0$.
If $d_G(v_0)+d_G(v_1)\geq n$, then $G$ is either pancyclic, or bipartite, or missing only a cycle of length $n-1$.
If $d_G(v_0)+d_G(v_1)> n$, then $G$ is pancyclic.
\end{theorem}
The following three lemmas capture obvious remarks.

\begin{lemma}\label{lemma1}
If $G$ is a graph of order $n$ with a cycle $C:v_0v_1\ldots v_{n-1}v_0$ of length $n-1$ 
and the vertex $u$ that does not lie on $C$ has degree more than $\frac{n-1}{2}$, then $G$ is pancyclic.
\end{lemma}
{\it Proof:} It follows easily by the pigeon-hole principle that, 
for every integer $\ell$ with $3\leq \ell\leq n$, 
there is some index $i$ such that the vertices $v_i$ and $v_{i+\ell-1}$ are both neighbors of $u$,
where we identify indices modulo $n-1$,
that is, $G$ contains the cycle $v_iv_{i+1}\ldots v_{i+\ell-1}uv_i$ of length $\ell$
$\Box$

\begin{lemma}\label{lemma2}
If $G$ is a graph with a hamiltonian cycle $v_0v_1\ldots v_{n-1}v_0$ and maximum degree $\Delta$, 
then $S(G)$ contains $\left\lceil\frac{\Delta}{2}\right\rceil$ cycle lengths between $\frac{n+2}{2}$ and $n$.
\end{lemma}
{\it Proof:} For every neighbor $v_i$ of $v_0$ that is distinct from $v_1$ and $v_{n-1}$, 
one of the two cycles $v_0v_iv_{i+1}\ldots v_{n-1}v_0$ and $v_0v_iv_{i-1}\ldots v_1v_0$
has length between $\frac{n+2}{2}$ and $n-1$.
Since every such cycle length arises at most twice,
there are $\left\lceil\frac{\Delta-2}{2}\right\rceil$ cycle lengths between $\frac{n+2}{2}$ and $n-1$.
Together with the hamiltonian cycle, this implies the desired statement. $\Box$

\begin{lemma}\label{lemma3}
If $G$ is a connected graph with no independent set of order $3$, 
then $G$ has a hamiltonian path.
\end{lemma}
{\it Proof:} Let $P:v_1\ldots v_\ell$ be a longest path in $G$.
For a contradiction, we assume that $\ell<n$.
Since $G$ is connected, there is a vertex $u$ not on $P$ that has a neighbor, say $v_i$, on $P$.
Since $P$ is a longest path, $u$ is not a neighbor of $v_1$ or $v_\ell$.
Since $G$ has no independent set of order $3$,
the vertices $v_1$ and $v_\ell$ are adjacent.
Now $uv_iv_{i+1}\ldots v_\ell v_1v_2\ldots v_{i-1}$ is a path in $G$ 
contradicting the choice of $P$. $\Box$

\medskip

\noindent We proceed to the proof of our result.

\medskip

\noindent {\it Proof of Theorem \ref{theorem1}:}
Let $G$ be a claw-free hamiltonian graph of order $n$ and maximum degree $\Delta$ with $\Delta\geq 24$.
Let $C:v_0v_1\ldots v_{n-1}v_0$ be a hamiltonian cycle of $G$.
Let $d_G(v_0)=\Delta$.

First we assume that $\left\lceil\frac{3}{2}\Delta\right\rceil\geq n$.
For a contradiction, we may assume that $G$ is not pancyclic.
Since $\Delta>\frac{n-1}{2}$, 
Lemma \ref{lemma1} implies that $G$ does not contain the cycle $v_1v_2\ldots v_{n-1}v_1$ of length $n-1$,
that is, $v_1$ and $v_{n-1}$ are not adjacent.
Since $G$ is claw-free, 
every neighbor of $v_0$ that is distinct from $v_1$ and $v_{n-1}$, 
is adjacent to $v_1$ or $v_{n-1}$.
Since $v_0$ is adjacent to $v_1$ and $v_{n-1}$,
this implies $d_G(v_1)+d_G(v_{n-1})\geq \Delta$
and hence $d_G(v_0)+\max\{ d_G(v_1),d_G(v_{n-1})\}\geq \Delta+\left\lceil\frac{\Delta}{2}\right\rceil\geq n$.
Since $G$ is claw-free and $\Delta\geq 3$, the graph $G$ is not bipartite.
By Theorem \ref{theoremhasc}, we obtain that $G$ misses a cycle of length $n-1$.
This implies that $v_2$ and $v_{n-2}$ are not neighbors of $v_0$, 
that is, $v_1$ and $v_{n-1}$ both have two neighbors each that are not neighbors of $v_0$.
This implies
$d_G(v_1)+d_G(v_{n-1})\geq \Delta+2$
and hence $d_G(v_0)+\max\{ d_G(v_1),d_G(v_{n-1})\}\geq \Delta+\left\lceil\frac{\Delta}{2}\right\rceil+1>n$.
Now Theorem \ref{theoremhasc} yields a contradiction.

From now on we assume that $\left\lceil\frac{3}{2}\Delta\right\rceil< n$.

Since $G$ is claw-free, the subgraph $G[N_G(v_0)]$ of $G$ that is induced by the neighborhood $N_G(v_0)$ of $v_0$ 
is either connected or the disjoint union of two cliques. We consider two cases accordingly.

\medskip

\noindent {\bf Case 1} {\it $G[N_G(v_0)]$ is the disjoint union of two cliques $A$ and $B$.}

\medskip

\noindent The larger of the two cliques $A$ and $B$ together with $v_0$ 
forms a complete graph of order at least $\left\lceil\frac{\Delta}{2}\right\rceil+1$,
which implies $\left\{ 3,\ldots,\left\lceil\frac{\Delta}{2}\right\rceil+1\right\}\subseteq S(G)$.

Let $v_1\in A$.

First, we assume that $v_1$ and $v_{n-1}$ are adjacent, that is, $v_{n-1}\in A$.
By Lemma \ref{lemma1}, we may assume $\Delta\leq \frac{n-1}{2}$ in this case.
The path $v_1\ldots v_{n-1}$ contains a subpath $P$ 
between a vertex in $A$ and a vertex in $B$ 
whose internal vertices do not belong to $N_G[v_0]$.
Let $P$ have $p$ internal vertices.
In view of the structure of $G$, 
there are cycles in $G$ that contain $v_0$, $P$, and some further vertices of $A\cup B$ of all lengths 
between $p+3$ and $p+\Delta+1$, that is,
$\left\{ p+3,\ldots,p+\Delta+1\right\}\subseteq S(G)$.
If $p+3>\left\lceil\frac{\Delta}{2}\right\rceil+1$, 
this already implies 
$s(G)\geq \left\lceil\frac{\Delta}{2}\right\rceil-1+\Delta-1=\left\lceil\frac{3}{2}\Delta\right\rceil-2$.
Hence, we may assume $p+3\leq\left\lceil\frac{\Delta}{2}\right\rceil+1$,
which implies $\left\{ 3,\ldots,\Delta\right\}\subseteq S(G)$.
Since $\Delta\leq \frac{n-1}{2}$,
Lemma \ref{lemma2} implies $s(G)\geq 
\Delta-2+\left\lceil\frac{\Delta}{2}\right\rceil=\left\lceil\frac{3}{2}\Delta\right\rceil-2$.

From now on we assume that $v_1$ and $v_{n-1}$ are not adjacent, that is, $v_{n-1}\in B$.

Removing from the hamiltonian cycle $C$ 
the vertex $v_0$ and each of its edges that lies completely in either $A$ or $B$,
and splitting the arising paths at vertices in $A\cup B$,
results in a collection of paths of order either $1$ or at least $3$
whose endvertices are in $A\cup B$ and whose internal vertices do not belong to $N_G[v_0]$.
Whenever some vertex $v$ in $A\cup B$ is the endvertex of two of these paths of order at least $3$, 
say $P:vx_1\ldots x_p$ and $Q:vy_1\ldots y_q$,
then the claw-freeness of $G$ implies that the two neighbors $x_1$ and $y_1$ of $v$ on $P$ and $Q$ are adjacent. 
In such a situation, 
we replace the two paths $P$ and $Q$ within our collection 
by the single path $x_p\ldots x_1y_1\ldots y_q$,
that is, we shortcut $P\cup Q$ by removing $v$.
Iteratively repeating this replacement as long as possible and removing the paths of order $1$,
we obtain a collection ${\cal P}=\{ P_1,\ldots,P_k\}$ of vertex-disjoint paths of order at least $3$
whose endvertices are in $A\cup B$ and whose internal vertices do not belong to $\{ v_0\}\cup A\cup B$
such that every vertex of $G$ that does not belong to $\{ v_0\}\cup A\cup B$
is an internal vertex of some path in ${\cal P}$.
Since the paths in ${\cal P}$ have distinct endvertices in $A\cup B$, we have $k\leq \frac{\Delta}{2}$.
If $p_i$ denotes the number of internal vertices of $P_i$,
then $P_i$ has order $p_i+2$ and $p_1+\cdots+p_k=n-\Delta+1$.
Since $v_1\in A$ and $v_{n-1}\in B$, 
it follows easily that ${\cal P}$ contains a path, say $P_1$, between a vertex in $A$ and a vertex in $B$.

In view of the structure of $G$, 
it follows that, 
for every integer $i$ with $1\leq i\leq k$,
there are cycles in $G$ that contain $v_0$, $P_1,\ldots,P_i$, and some further vertices of $A\cup B$ of all lengths 
between $1+(p_1+2)+\cdots+(p_i+2)=1+2i+p_1+\cdots+p_i$ 
and $1+\Delta+p_1+\cdots+p_i$.
If the union of the following subsets of $S(G)$ 
\begin{eqnarray*}
&& \left\{ 3,\ldots,\left\lceil\frac{\Delta}{2}\right\rceil+1\right\}\\
&& \left\{ 3+p_1,\ldots,1+\Delta+p_1\right\}\\
&& \left\{ 5+p_1+p_2,\ldots,1+\Delta+p_1+p_2\right\}\\
&& \ldots\\
&& \left\{ 1+2k+p_1+\cdots+p_k,\ldots,1+\Delta+p_1+\cdots+p_k\right\}
\end{eqnarray*}
is a set of consecutive integers, 
then $1+\Delta+p_1+\cdots+p_k=n$ implies that $G$ is pancyclic.
If 
$\left\{ 3,\ldots,\left\lceil\frac{\Delta}{2}\right\rceil+1\right\}$
and 
$\left\{ 3+p_1,\ldots,1+\Delta+p_1\right\}$
are disjoint subsets of $S(G)$, then 
$s(G)\geq \left\lceil\frac{\Delta}{2}\right\rceil-1+\Delta-1=\left\lceil\frac{3}{2}\Delta\right\rceil-2$.
Hence we may assume that there is some index $i$ with $2\leq i\leq k$ such that 
the two sets $\left\{ 3,\ldots,1+\Delta+p_1+\cdots+p_{i-1}\right\}$
and
$\left\{ 1+2i+p_1+\cdots+p_i,\ldots,1+\Delta+p_1+\cdots+p_i\right\}$
are disjoint subsets of $S(G)$.
This implies
\begin{eqnarray*}
s(G) & \geq  & (\Delta+p_1+\cdots+p_{i-1}-1)+(\Delta-2i+1)\\
&\geq & 2\Delta-i-1\\
& \geq & 2\Delta-k-1\\
& \geq & \frac{3}{2}\Delta-1\\
& > & \left\lceil\frac{3}{2}\Delta\right\rceil-2,
\end{eqnarray*}
which concludes the proof in this case.

\medskip

\noindent {\bf Case 2} {\it $G[N_G(v_0)]$ is connected.}

\medskip

\noindent Let 
\begin{eqnarray*}
A & = & \{ v_1\}\cup ((N_G(v_0)\cap N_G(v_1))\setminus N_G(v_{n-1})),\\
B & = & \{ v_{n-1}\}\cup ((N_G(v_0)\cap N_G(v_{n-1}))\setminus N_G(v_1)),\mbox{ and }\\
X & = & N_G(v_0)\cap N_G(v_1)\cap N_G(v_{n-1}).
\end{eqnarray*}
Since $G$ is claw-free, Lemma \ref{lemma3} implies that $G[N_G(v_0)]$ has a hamiltonian path, 
which implies $\left\{ 3,\ldots,\Delta+1\right\}\subseteq S(G)$.
If $\Delta+1\leq \frac{n+2}{2}$, then Lemma \ref{lemma2} implies
$s(G)\geq \Delta-1+\left\lceil\frac{\Delta}{2}\right\rceil-1=\left\lceil\frac{3}{2}\Delta\right\rceil-2$.
Hence we may assume that $\Delta\geq \frac{n+1}{2}$.
By Lemma \ref{lemma1}, we may assume that the vertices $v_1$ and $v_{n-1}$ are not adjacent.
Similarly, if we have $v_i\in B\cup X$ and $v_{i+1}\in A\cup X$ for some index $i$, 
then $v_1v_2\ldots v_iv_{n-1}v_{n-2}\ldots v_{i+1}v_1$ is a cycle of length $n-1$ in $G$
and $d_G(v_0)\geq \frac{n+1}{2}>\frac{n-1}{2}$ implies that $G$ is pancyclic.
Therefore, for every index $i$, 
we may assume that $v_i\not\in B\cup X$ or $v_{i+1}\not\in A\cup X$. 

Clearly, we may assume that $G$ is not pancyclic and 
consider the positive integer $\ell$ chosen in such a way that 
$n-\ell$ is the largest element less than $n$ 
that is missing from $S(G)$.

We orient all edges of the hamiltonian cycle $C$ 
according to the order $v_0,v_1,\ldots ,v_{n-1},v_0$.
This results in an oriented cycle $\vec{C}$.
Removing from $\vec{C}$ 
the vertex $v_0$ and each of its oriented edges that lies completely in either $A$, or $B$, or $X$,
and splitting the arising paths at vertices in $A\cup B\cup X$,
results in a collection $\vec{{\cal P}}$ of oriented paths
whose endvertices are in $A\cup B\cup X$ 
and whose internal vertices do not belong to $N_G[v_0]$.
We refer to paths in $\vec{{\cal P}}$ 
that are between two distinct of the sets $A$, $B$, and $X$ as {\it crossing} 
and to paths in $\vec{{\cal P}}$ 
from $A$ to $A$, or from $B$ to $B$, or from $X$ to $X$ as {\it non-crossing}.
Unlike in Case 1, the crossing paths may have order $2$.
By construction, a non-crossing path has at least one internal vertex.

Since $G$ is claw-free and $A$ and $B$ are non-empty by definition,
the sets $A$ and $B$ are cliques in $G$.

The above observation,
that is, $v_i\not\in B\cup X$ or $v_{i+1}\not\in A\cup X$ for every index $i$,
implies that no edge of $C$ lies completely in $X$, 
that is, in every vertex in $X$, one of the oriented paths in $\vec{{\cal P}}$ ends
and another one begins.

An oriented path in $\vec{{\cal P}}$ that begins in $B\cup X$ and ends in $A\cup X$ is called {\it reversed}.
Note that a path in $\vec{{\cal P}}$ from $X$ to $X$ is non-crossing as well as reversed.
If $\vec{P}:v_iv_{i+1}\ldots v_{i+p+1}$ is a reversed path with $p$ internal vertices, 
then the graph $G'=G-\{ v_{i+1},\ldots,v_{i+p}\}$ of order $n-p$
has a cycle $v_1v_2\ldots v_iv_{n-1}v_{n-2}\ldots v_{i+p+1}v_1$ of length $(n-p)-1$.
Since $d_{G'}(v_0)=d_G(v_0)\geq \frac{n+1}{2}$,
Lemma \ref{lemma2} and the choice of $\ell$ imply that 
$\{ 3,\ldots,n-p\}\cup \{ n-\ell+1,\ldots,n\}\subseteq S(G)$.
By the choice of $\ell$, we have $n-p<n-\ell$ and thus $s(G)\geq (n-p-2)+\ell$.
Clearly, we may assume $s(G)<\left\lceil\frac{3}{2}\Delta\right\rceil-2$ 
and thus $p\geq n-\left\lceil\frac{3}{2}\Delta\right\rceil+\ell+1$.

Let $G$ have $r$ reversed paths.

First we assume that $r\geq 1$.

We mark the vertices in $N_G(v_0)$ as follows:
\begin{itemize}
\item Every vertex in $A\cup B$ that lies on a crossing path in $\vec{{\cal P}}$
is marked with ``$M_1$''.
\item Every vertex in $X$ is also marked with ``$M_1$''.
\item Every vertex in $A\cup B$ that lies on a non-crossing path in $\vec{{\cal P}}$
that has at least two internal vertices, is marked with ``$M_2$''.
\end{itemize}
Note that vertices in $A\cup B$ may receive both marks $M_1$ and $M_2$.
Let $m_1$ and $m_2$ denote the number of vertices marked with $M_1$ and $M_2$, respectively.

A vertex that does receive none of the two marks is {\it free}.
Note that a vertex $u$ is free if it lies in $A\cup B$ 
and every path in $\vec{{\cal P}}$
to which $u$ belongs, 
is non-crossing and has exactly one internal vertex.

Let $G$ have $f$ free vertices.

Note that for every free vertex $v_j$ that belongs to $D\in \{ A,B\}$, 
the hamiltonian cycle $C$ contains a subpath $P:v_i\ldots v_k$ such that 
\begin{itemize}
\item $i<j<k$, that is, $v_j$ lies on $P$,
\item $v_i$ and $v_k$ are non-free vertices in $D$, and
\item every internal vertex of $P$
\begin{itemize}
\item is either a free vertex from $D$ 
\item or is the unique internal vertex of a path in $\vec{{\cal P}}$ of length $2$. 
\end{itemize}
\end{itemize}
Note that the vertices in $V(P)\cap D$ form a clique.

This structure around the $f$ free vertices of $G$ easily implies that 
$\{ n-f,\ldots,n\}\subseteq S(G)$ and,
by the choice of $\ell$, it follows that $f<\ell$.

By the marking process, 
the number of internal vertices 
of non-crossing paths in $\vec{{\cal P}}$ 
with at least $2$ internal vertices
that are not between $X$ and $X$, 
is at least $m_2$.
Since every of the $r$ reversed paths of $G$ has $n-\left\lceil\frac{3}{2}\Delta\right\rceil+\ell+1$ internal vertices,
we obtain 
\begin{eqnarray}\label{e1}
m_2+\left(n-\left\lceil\frac{3}{2}\Delta\right\rceil+\ell+1\right)r\leq n-\Delta-1.
\end{eqnarray}
Removing from $C$ the vertex $v_0$ and the internal vertices of the $r$ reversed paths,
results in a collection ${\cal Q}$ of $r+1$ paths, 
one of which contains $v_1$ and one of which contains $v_{n-1}$.
The path in ${\cal Q}$ that contains $v_1$ contains at most $4$ vertices marked with $M_1$;
one from $A$, one from $X$, and two from $B$.
Similarly, the path in ${\cal Q}$ that contains $v_{n-1}$ contains at most $4$ vertices marked with $M_1$.
All remaining paths in ${\cal Q}$ contain at most $5$ vertices marked with $M_1$;
two from $A$, one from $X$, and two from $B$.
Altogether, there are at most $5r+3$ vertices marked with $M_1$.
This implies that there are at least $\Delta-(5r+3)-m_2$ free vertices, that is,
\begin{eqnarray}\label{e2}
\ell-1 \geq f\geq \Delta-(5r+3)-m_2.
\end{eqnarray}
Using $n-\left\lceil\frac{3}{2}\Delta\right\rceil\geq 1$ and $\ell\geq 1$, 
(\ref{e1}) implies 
\begin{eqnarray}
m_2 & \leq & n-\Delta-1-\left(n-\left\lceil\frac{3}{2}\Delta\right\rceil+\ell+1\right)r\nonumber\\
& = & \left\lceil\frac{\Delta}{2}\right\rceil-\ell-2-\left(n-\left\lceil\frac{3}{2}\Delta\right\rceil+\ell+1\right)(r-1)\label{e3}\\
& \leq & \left\lceil\frac{\Delta}{2}\right\rceil-3-3(r-1).\nonumber
\end{eqnarray}
Using $m_2\geq 0$, this implies $r\leq \frac{1}{3}\left\lceil\frac{\Delta}{2}\right\rceil$.
Combining (\ref{e2}) and (\ref{e3}) and using $n-\left\lceil\frac{3}{2}\Delta\right\rceil\geq 1$ and $\ell\geq 1$, we obtain 
\begin{eqnarray*}
\ell-1 & \geq &
\Delta-(5r+3)-\left(\left\lceil\frac{\Delta}{2}\right\rceil-\ell-2-\left(n-\left\lceil\frac{3}{2}\Delta\right\rceil+\ell+1\right)(r-1)\right)\\
& = &
\left\lfloor\frac{\Delta}{2}\right\rfloor+(\ell-1)-5+\left(n-\left\lceil\frac{3}{2}\Delta\right\rceil+\ell-4\right)(r-1)\\
& \geq &
\left\lfloor\frac{\Delta}{2}\right\rfloor+(\ell-1)-5-2(r-1),
\end{eqnarray*}
which implies
$r\geq \frac{1}{2}\left(\left\lfloor\frac{\Delta}{2}\right\rfloor-3\right)$.
Since 
$\frac{1}{2}\left(\left\lfloor\frac{\Delta}{2}\right\rfloor-3\right)\leq \frac{1}{3}\left\lceil\frac{\Delta}{2}\right\rceil$
implies $\Delta\leq 23$, which completes the proof in the case $r\geq 1$. 

Next we assume that $r=0$.
By the definition of a reversed path, this implies that $X$ contains at most one vertex and that 
the hamiltonian cycle $C$, 
traversed in the order $v_0,v_1,v_2,\ldots$ and so on, 
first visits all vertices of $A$,
then all vertices of $X$, 
and then all vertices of $B$.
Therefore, the cycle $C$ contains exactly one subpath from a vertex in $A$, say $v_a$, to a vertex in $B$, say $v_b$,
and, if $X$ is non-empty, then the unique vertex in $X$ is an internal vertex of this path.

Manipulating $C$ similarly as in Case 1, 
we obtain the existence of two sets
$A'$ with $A'\in \{ A,A\setminus \{ v_a\}\}$
and 
$B'$ with $B'\in \{ B,B\setminus \{ v_b\}\}$,
and a collection ${\cal R}=\{ R_1,\ldots,R_t\}$ of vertex-disjoint paths 
whose endvertices are in $A'\cup B'$ and whose internal vertices do not belong to $\{ v_0\}\cup A'\cup B'$
such that every vertex of $G$ that does not belong to $\{ v_0\}\cup A'\cup B'$
is an internal vertex of some path in ${\cal R}$.
Note that the vertex $v_a$ may play a special role if it is adjacent in $C$ to $x$ or $v_b$
and its second neighbor on $C$ does not belong to $A$.
In this situation we cannot necessarily perform a shortcut as in Case 1.
Instead we may have to remove $v_a$ from $A$ to obtain $A'$ 
and consider $v_a$ an internal vertex of some path in ${\cal R}$.
Similar comments apply to $v_b$.

Since $C$ contains exactly one subpath from $A$ to $B$,
the collection ${\cal R}$ contains exactly one path, say $R_1$, between $A'$ and $B'$.
Let $q_i$ denote the number of internal vertices of $R_i$.
The vertex in $X$, if it exists, is an internal vertex of $R_1$.
As explained above, $v_a$ and $v_b$ may be internal vertices of $R_1$,
and, if $v\in \{ v_a,v_b\}$ is an internal vertex of $R_1$,
then $R_1$ contains a neighbor of $v$ that does not belong to $N_G[v_0]$ as an internal vertex.
If $\Delta'$ denotes $|A'\cup B'|$, then we obtain $\Delta-3\leq \Delta'\leq \Delta$.
Furthermore, 
if $\Delta'=\Delta-3$, then $q_1\geq 5$,
if $\Delta'=\Delta-2$, then $q_1\geq 3$, 
if $\Delta'=\Delta-1$, then $q_1\geq 1$, and
if $\Delta'=\Delta$, then $q_1\geq 0$.
This implies that $\frac{3}{2}\Delta'+q_1\geq \frac{3}{2}\Delta$ 
unless $\Delta'=\Delta-1$ and $q_1=1$,
which implies $A'=A$, $B'=B$, and $|X|=1$.
For $i\geq 2$, we have $q_i\geq 1$.
Since the paths in ${\cal R}$ have distinct endvertices in $A'\cup B'$, we have $t\leq \frac{\Delta'}{2}$.

In view of the structure of $G$, 
it follows that, 
for every integer $i$ with $1\leq i\leq t$,
there are cycles in $G$ that contain $v_0$, $R_1,\ldots,R_i$, and some further vertices of $A'\cup B'$ of all lengths 
between $1+(q_1+2)+\cdots+(q_i+2)=1+2i+q_1+\cdots+q_i$ 
and $1+\Delta'+q_1+\cdots+q_i$.
If the union of the following subsets of $S(G)$ 
\begin{eqnarray*}
&& \left\{ 3,\ldots,\Delta+1\right\}\\
&& \left\{ 3+q_1,\ldots,1+\Delta'+q_1\right\}\\
&& \left\{ 5+q_1+q_2,\ldots,1+\Delta'+q_1+q_2\right\}\\
&& \ldots\\
&& \left\{ 1+2t+q_1+\cdots+q_t,\ldots,1+\Delta'+q_1+\cdots+q_t\right\}
\end{eqnarray*}
is a set of consecutive integers, 
then $1+\Delta'+q_1+\cdots+q_t=n$ implies that $G$ is pancyclic.
If 
$\left\{ 3,\ldots,\Delta+1\right\}$
and 
$\left\{ 3+q_1,\ldots,1+\Delta'+q_1\right\}$
are disjoint subsets of $S(G)$, then 
$s(G)\geq \Delta-1+\Delta'-1\geq 2\Delta-5\geq \left\lceil\frac{3}{2}\Delta\right\rceil-2$.
Hence we may assume that there is some index $i$ with $2\leq i\leq t$ such that 
the two sets $\left\{ 3,\ldots,1+\Delta'+q_1+\cdots+q_{i-1}\right\}$
and
$\left\{ 1+2i+q_1+\cdots+q_i,\ldots,1+\Delta'+q_1+\cdots+q_i\right\}$
are disjoint subsets of $S(G)$.
This implies
\begin{eqnarray*}
s(G) &\geq & (\Delta'+q_1+\cdots+q_{i-1}-1)+(\Delta'-2i+1)\\
& = & 2\Delta'+q_1+(q_2+\cdots+q_{i-1})-2i\\
& \geq & 2\Delta'+q_1+(i-2)-2i\\
& = & 2\Delta'+q_1-i-2\\
& \geq & 2\Delta'+q_1-t-2\\
& \geq & 2\Delta'+q_1-\frac{\Delta'}{2}-2\\
& = & \frac{3}{2}\Delta'+q_1-2.
\end{eqnarray*}
This already implies the desired bound unless
$\Delta'=\Delta-1$,
$q_1=1$,
$A'=A$, 
$B'=B$, 
$|X|=1$,
$q_2=\ldots=q_{t-1}=1$, and 
$i=t=\frac{\Delta'}{2}$.
Since $t\geq 3$ and the order of the paths $R_2,\ldots,R_t$ was arbitrary,
we obtain, by symmetry, that $q_t=1$ as well.
Altogether, every vertex in $A\cup B$ is the endvertex of exactly one of the paths in ${\cal R}$
and all these paths have exactly $1$ internal vertex.
This implies that 
$n=|\{ v_0\}|+|A\cup B|+\frac{|A\cup B|}{2}=1+\frac{3}{2}(\Delta-1)=\frac{3}{2}\Delta-\frac{1}{2}$,
which contradicts $\left\lceil\frac{3}{2}\Delta\right\rceil< n$
and completes the proof. $\Box$

\medskip

\noindent The condition $\Delta\geq 24$ in the statement of Theorem \ref{theorem1} is an artefact of our proof;
we chose to avoid the kind of tedious detailed case analysis, 
which would easily allow to reduce this bound.
In fact, our proof already implies that graphs, 
which do not satisfy the conclusion of Theorem \ref{theorem1},
satisfy $\frac{n+1}{2}\leq \Delta\leq 23$,
that is, there are at most finitely many such graphs.

\pagebreak

\end{document}